\documentclass[12pt]{amsart} 
\usepackage{lmodern}
\usepackage{amsmath}
\usepackage{amssymb}
\usepackage{latexsym}
\def\Mb{{\mathbf M}}

\def\Er{{\mathbb E}}

\def\Pr{{\mathbb P}}

\def\Rr{{\mathbb R}}

\def\Ac{{\mathcal{A}}}

\def\Dc{{\mathcal{D}}}
\def\Ec{{\mathcal{E}}}
\def\Fc{{\mathcal{F}}}

\def\Rc{{\mathcal{R}}}

\def\Tc{{\mathcal{T}}}

\def\one{{\rm \bf 1}}
\def\({\left(}     
\def\){\right)}    
\def\[{\left[}     
\def\]{\right]}

\newtheorem{definition}{Definition}
\newtheorem{theorem}{Theorem}
\newtheorem{lemma}{Lemma}
\newtheorem{proposition}{Proposition}
\newtheorem{remark}{Remark}

\begin{document}
\author{Freddy Delbaen}
\address{Departement f\"ur Mathematik, ETH Z\"urich, R\"{a}mistrasse
101, 8092 Z\"{u}rich, Switzerland, also Institut f\"ur Mathematik,
Universit\"at Z\"urich, Winterthurerstrasse 190,
8057 Z\"urich, Switzerland}
\email{delbaen@math.ethz.ch}
\author{Chitro Majumdar}
\address{RSRL, Chitro Majumdar, Jumeirah Beach Residence (JBR), Dubai, United Arab Emirates}
\email{chitro.majumdar@rsquarerisklab.com; \newline chitromajumdar@icloud.com}
\title{Approximation with Independent Variables\footnote{AMS-classification primary 60G35, secondary 68T99, 93E11, 94A99}}
\maketitle
\abstract  Given a square integrable m-dimensional random variable $X$ on a probability space $(\Omega.\Fc,\Pr)$ and a sub sigma algebra $\Ac$, we show that there exists another m-dimensional random variable $Y$, independent of $\Ac$ and minimising the $L^2$ distance to $X$.  Such results have an importance to fairness and bias reduction in Artificial Intelligence, Machine Learning and Network Theory. 
\endabstract

 \section{Notation and Preliminaries}
We use standard probabilistic notation, $(\Omega,\Fc,\Pr)$ is a probability space  and $\Ac \subset\Fc$ is a sub $\sigma-$algebra. All random variables will be square integrable and all norms will be the $L^2$ norm.  For random variables $X\colon \Omega \rightarrow \Rr^m$, the norm is then defined as $\Vert X\Vert^2=\int |X|^2\Pr[d\omega]$, where for $x\in \Rr^m$, $|x|$ denotes the Euclidean norm of $x$.  The problem we want to solve is the following. Given a random variable  $X\colon\Omega \rightarrow \Rr^m$,  $X \in L^2 $, can we find a random variable $Y\colon \Omega\rightarrow \Rr^m$ which is independent of  the $\sigma-$algebra $\Ac$ and which minimises $\Vert X-Y\Vert$ among all $m-$dimensional random variables that are independent of $\Ac$.

The motivation to study this problem found its origin in artificial intelligence and neural networks. If a random outcome $X$ is observed that is subject to perturbations described by the information structure or sigma algebra $\Ac$, then we would like to replace $X$ by a random variable $Y$ that is independent of these perturbations or bias. We use the $L^2$ distance as a measure for approximation but the method of proof of the existence of $Y$ also works for $L^p, 1\le p < \infty$ distances. However in this case the relation with correlation coefficients disappears.

The idea of the proof is to use results of transport theory and the Wasserstein metric.  We do not need the deeper results of this theory so the reader can look up the necessary tools in any book on this topic, see for instance \cite{Thorpe}. To find an independent random variable that solves our problem we will use a parametrised form of the solution of suitable Monge-Kantorovitch problems. To solve the necessary measurability issues we will state and prove some (probably known) measurability results on the dependence of optimal transport plans.  In case the dimension $m=1$, the use of transport theory can be simplified and can be replaced -- as is well known from transport theory -- by the use of commonotonicity. The optimisation of the $L^2$ distance is then based on the inequalities 368 page 261 and 378 page 278 in \cite{HLP}.
 
To avoid possible complications with sets of measure zero,  we will suppose that ${\Ac}$ contains all null sets of $\Fc$. Since we are looking for a random variable that is independent of $\Ac$ we suppose that there is  a $U\colon \Omega \rightarrow [0,1]$ uniformly distributed and independent of $\Ac$.  This condition is equivalent to the condition that $\Fc$ is atomless with respect to $\Ac$, see \cite{del-atom}.  It guarantees the existence of random variables with arbitrary distribution.

We now recall some facts about conditional distributions and independence.  The reader can look up the details in \cite{Rao}, page 126, theorem 5 where the theorems are proved for $\Rr^m$ (the cases we will use), but they remain valid for Polish spaces, see \cite{Bill}.  
\begin{lemma}Let $E$ be a Polish space equipped with its Borel sigma-algebra $\Ec$. Let $Z\colon \Omega\rightarrow E$ be a measurable function.  Furthermore let $\Ac\subset \Fc$ be a  sub sigma-algebra.  There exists a kernel $\mu\colon \Omega\times \Ec$ such that 
\begin{enumerate}
\item for each $F\in\Ec$, the mapping $\omega\rightarrow \mu(\omega,F)$ is $\Ac$ measurable
\item for almost every $\omega\in \Omega$ the mapping $F\rightarrow \mu(\omega,F)$ defines a probability on $\Ec$
\item  for each bounded measurable function $h\colon E\rightarrow \Rr$, the conditional expectation is given by  $\Er[h(Z)\mid\Ac](\omega)=\int_\Ec h(x)\mu(\omega,dx)$ almost everywhere
\end{enumerate}
\end{lemma}
\begin{lemma} With the notation of the previous lemma we have that $Z$ is independent of  $\Ac$  if and only if almost everywhere $\mu(\omega,.)=\nu$, where $\nu=\Pr\circ Z^{-1}$ the distribution of $Z$.
\end{lemma}
\begin{lemma} $Z$ is independent of  $\Ac$ if and  only if for each $A\in \Ac$  and each bounded continuous function $h\colon \Rr\rightarrow \Rr$ we have $\Er[h(Z)\one_A]=\Pr[A]\,\Er[h(Z)]$.  For sets the situation is easier,  a set $B\in \Fc$ is independent of $\Ac$ if and only  if $\Er[\one_B\mid\Ac]=\Pr[B]$.\end{lemma}

 \section{Realisation of measures with random variables}
 
It is well known that for a Polish space $(E,\Ec)$ an arbitrary probability measure $\alpha$ on $\Ec$ can be realised with a random variable, a Borel measurable function $\psi\colon [0,1]\rightarrow E$, meaning that the distribution of $\psi$ equals $\alpha$. We need a parametrised version of this result.  The following theorem, due to Blackwell and Dubins, \cite{Black-Dub} says a lot more than we need.
 \begin{theorem}  For a Polish space $(E,\Ec)$, let $\Mb$ be the set of probability measures on $\Ec$. There exists a Borel measurable function $h\colon \Mb\times[0,1]\rightarrow E$ such that for each $\alpha\in \Mb$, the distribution of $h(\alpha,.)$ ($[0,1]$ is equipped with the Lebesgue measure $\lambda$) is precisely $\alpha$.  Furthermore if $\alpha_n\rightarrow\alpha$ is a weak$^*$ converging sequence then $h(\alpha_n,t)\rightarrow h(\alpha,t)$ $\lambda$ almost surely. 
 \end{theorem}\label{h}
 \begin{remark}  {\rm We do not need the continuity part of the theorem.  We only need the existence and the Borel measurability of $h$.}
\end{remark}
\begin{remark} {\rm   The Blackwell-Dubins theorem is a refinement of the Skorohod theorem which says that if $\mu_n,n\ge 1$ is a sequence of probability measures on a Polish space $E$, converging weak$^*$ to a probability measure $\mu$ then there exists a sequence of Borel measurable functions $X_n,n\ge 1$ and a Borel measurable function $X$ such that $X_n\colon [0,1]\rightarrow E$, $X_n\rightarrow X$ almost everywhere and the distribution of $X_n$ is $\mu_n$.  Of course the distribution of $X$ is then $\mu$.  This allows to replace the weak$^*$ convergence by a pointwise convergence of measurable functions.  We will, without mentioning, use this trick  when we use the monotone convergence theorem, Fatou's lemma and truncation arguments for weak$^*$ convergence.}
\end{remark} 
 \section{The Wasserstein Metric and Transport Problems}
 
 For two probability measures $\mu,\nu$ on $\Rr^m$, each having second moments i.e. $\int_{\Rr^m}|x|^2 d\mu + \int_{\Rr^m}|x|^2 d\nu<\infty$, we denote by $\Pi(\mu,\nu)$ the set of probability measures $\gamma$ on $\Rr^m\times \Rr^m$ such that the marginals of $\gamma$ are precisely $\mu$ and $\nu$.  Obviously this set is closed in the weak$^*$ topology of probability measures on $\Rr^m\times\Rr^m$. The set $K(\mu,\nu)$ denotes the subset of $\Pi(\mu,\nu)$ realising the minimum of $\int_{\Rr^m\times \Rr^m}|x-y|^2 d\gamma$ on $\Pi(\mu,\nu)$.  The set $K(\mu,\nu)$ was shown to exist by Kantorovitch and is the object of transport theory, \cite{Thorpe}.  To make the statements easier to read, let us denote by $\Mb_2(k)$ the set of probability measures on $\Rr^k$ having finite second moment. There are two topologies on $\Mb_2(k)$ that are of interest to us.  The first one is the weak$^*$ topology, the second one is the topology asking the convergence of $\int f d\gamma$ where $f$ runs through the set of functions that are quadratically bounded, i.e. there is a constant $C$ such that for all $x\in\Rr^k$ we have $|f(x)|\le C(1+|x|^2)$.  This topology which we denote by $\Tc_2$, turns $\Mb_2(k)$ into a Polish space. As easily seen $\mu_n\rightarrow \mu$ for $\Tc_2$ if and only if $\mu_n\rightarrow\mu$ weak$^*$ and $\int_{\Rr^k}|x|^2d\mu_n\rightarrow \int_{\Rr^k}|x|^2d\mu$. The square of the Wasserstein distance $W(\mu,\nu)$ between two elements $\mu,\nu\in\Mb_2(m)$ is defined as $W^2(\mu,\nu)=\int_{\Rr^m\times\Rr^m}|x-y|^2d\gamma$ where $\gamma\in K(\mu,\nu)$.
 \begin{theorem} With the notation introduced in the preceding paragraph, the set
 $$
 \left\{(\mu,\nu,\gamma)\mid \mu,\nu\in \Mb_2(m),\gamma\in K(\mu,\nu)   \right\}
 $$
 is closed in $\Mb_2(m)\times\Mb_2(m)\times\Mb_2(2m)$, equipped with the product topology of the corresponding $\Tc_2$ topologies.  For each $(\mu,\nu)\in \Mb_2(m)\times\Mb_2(m)$, the set $K(\mu,\nu)$ is compact for $\Tc_2$. The Wasserstein distance is continuous on $\Mb_2(m)\times\Mb_2(m)$ for the product topology of the $\Tc_2$ topologies.
 \end{theorem}
 {\bf Proof} These results are standard but for completeness we will give a proof, \cite{Thorpe}. We will prove the closedness of the graph in such a way that it includes a proof of compactness.   Let $\mu_n\rightarrow \mu, \nu_n\rightarrow \nu$ for the $\Tc_2$ topology. Let $\gamma_n$ be a sequence with $\gamma_n\in K(\mu_n,\nu_n)$. We will show that the sequence contains a subsequence that converges to an element of $K(\mu,\nu)$. That will prove closedness of the graph and compactness of $K(\mu,\nu)$.  We have $\int |(x,y)|^2d\gamma_n= \int |x|^2d\mu_n+\int |y|^2d\nu_n\rightarrow \int |x|^2d\mu+\int |y|^2d\nu<\infty$. It follows that the sequence $(\gamma_n)_n$ is relatively weak$^*$ compact.  Without loss of generality we may suppose that $\gamma_n\rightarrow \gamma_0$ for the weak$^*$ topology. Obviously $\gamma_0\in \Pi(\mu,\nu)$. Using truncation and the monotone convergence theorem we have $\int |(x,0)|^2 \,d\gamma_0=\int |x|^2\,d\mu,
 \int |(0,y)|^2 \,d\gamma_0=\int |y|^2\,d\nu$. This shows that $\gamma_n\rightarrow \gamma_0$ for the $\Tc_2$ topology. The dominated convergence theorem shows that $\int|x-y|^2\,d\gamma_n\rightarrow \int |x-y|^2\,d\gamma_0$. In other words $\lim_n W^2(\mu_n,\nu_n)=\int |x-y|^2\,d\gamma_0\ge W^2(\mu,\nu)$.  But the triangle inequality for the Wasserstein distance as well as its relation to weak$^*$ convergence shows that $W(\mu_n,\nu_n)\le W(\mu_n,\mu)+W(\mu,\nu)+W(\nu,\nu_n)$ and by taking limits for $n\rightarrow \infty$ we get $\limsup W(\mu_n,\nu_n)\le W(\mu,\nu)$.  This in turn gives $\gamma_0\in K(\mu,\nu)$.\hfill$\Box$
 
\begin{remark}{\rm  The set $\Mb_2(k)$ equipped with the weak$^*$ topology is a Lusin space (because $\Tc_2$ is a finer topology).  The Borel sets for $\Tc_2$ and for the weak$^*$ topology are therefore the same, see \cite{Par} and hence the set $$ \left\{(\mu,\nu,\gamma)\mid \mu,\nu\in \Mb_2(m),\gamma\in K(\mu,\nu)   \right\}$$ is a Borel set when the product space is equipped with the product of the weak$^*$ topologies.  An element $\gamma\in K(\mu,\nu)$ is sometimes called a transport plan.}\end{remark}
The following theorem is a direct application of Novikov's selection theorem for Borel sets with compact sections, \cite{Sri}.
\begin{theorem}  There is a Borel measurable map
$$
\kappa\colon \Mb_2(m)\times\Mb_2(m)\rightarrow \Mb_2(2m)
$$
such that for every $\mu,\nu\in\Mb_2(m)$, $\kappa(\mu,\nu)\in K(\mu,\nu)$.
\end{theorem}\label{kappa}
\begin{remark}{\rm It is easy to find examples showing that a continuous selection does not exist.}
\end{remark}
\section{Obvious Properties of the solution}

\begin{lemma} Let $\Ac\subset\Fc$  be a sub sigma-algebra and let $X\in L^2$.  If $Y$ is independent of $\Ac$ and  satisfies
$$
\Vert  X-Y\Vert=\inf\{ \Vert X-Z\Vert\mid Z\text{ independent of }\Ac\},
$$
then $\Er[X]=\Er[Y]$.
\end{lemma}
{\bf Proof }  This is rather obvious since $\Rr^m\rightarrow\Rr_+, \alpha\rightarrow\Vert X-Y-\alpha\Vert^2$ attains its minimum for $\alpha=\Er[X-Y]$.\hfill$\Box$
\begin{lemma} Let $\Ac\subset\Fc$  be a sub sigma-algebra and let $X\in L^2$.  If $X$ is $\Ac-$measurable then
$$
\Vert  X-\Er[X]\Vert=\inf\{ \Vert X-Z\Vert\mid Z\text{ independent of }\Ac\}.
$$
In other words the solution of the optimisation problem (as in the previous lemma) is $Y=\Er[X]$.
\end{lemma}
{\bf Proof} Because  of the previous lemma we can change the optimisation in the following way
$$\inf\{ \Vert X-Z\Vert\mid Z\text{ independent of }\Ac\text{ and }\Er[Z]=\Er[X]\}.$$
We then find
$$
\Vert X-Z\Vert^2= \Vert X-\Er[X]\Vert^2-2\langle X-\Er[X],Z-\Er[Z]\rangle + \Vert Z-\Er[Z]\Vert^2.
$$
The middle term is zero because of independence of $X$ and $Z$. The optimal choice for $Z$ is therefore  $Z=\Er[Z]=\Er[X]$.\hfill$\Box$

The statement in the following theorem can be void since we do not know if an optimal solution exists.  However in the next section  we will show by construction that an optimal solution exists under some additional hypothesis.
\begin{theorem} The optimal choice in 
$$\inf\{ \Vert X-Z\Vert\mid Z\text{ independent of }\Ac\}$$
is given by $Y=V+\Er[X]$ where $V$ is the optimal choice for the random variable $X-\Er[X\mid \Ac]$.
\end{theorem}
{\bf Proof } Take $V\in L^2$ and independent of $\Ac$. Let us write $Y=V+\Er[X]$ where as we know $\Er[X]=\Er[Y]$ for the relevant random variables. Then we write
$$
X-Y=\(X-\Er[X\mid\Ac] - V\) + (\Er[X\mid\Ac]-\Er[X]),
$$
which is a sum of two orthogonal variables.  Indeed $\(X-\Er[X\mid\Ac]\)\perp L^2(\Ac)$ and $V\perp L^2_0(\Ac)$, the space of $m-$dimensional $\Ac-$measurable random variables with expctation equal to zero. Hence
$$
\Vert X-Y\Vert^2=\Vert X-\Er[X\mid\Ac] - V\Vert^2 +\Vert \Er[X\mid\Ac]-\Er[X]\Vert^2.
$$
Therefore $\Vert X-Y\Vert^2$ is minimal if $\Vert X-\Er[X\mid\Ac] - V\Vert^2$ is minimal.
\hfill$\Box$
\begin{remark}{\rm We know that being uncorrelated is different from being independent.  In the preceding theorem the random variable $V$ is the best independent approximation  of the  random variable $X-\Er[X\mid\Ac]$  which is uncorrelated to the variables in $L^2(\Ac)$.  So the initial problem is in  fact equivalent to  the replacement  of uncorrelated random variables  by independent random variables.
}\end{remark}

\section{The Construction of the Solution}
 The idea is to ``decompose" $X$ along the atoms of $\Ac$.  This is done by constructing a factorisation of $X$ through a product space where the first factor  is $\Omega$ equipped with the sigma-algebra $\Ac$.  Then on each atom of $\Ac$ we will find $Y$ as the solution of a Monge-Kantorovitch problem.

 \begin{enumerate}
 \item Let $\Phi\colon \Omega \rightarrow\Omega\times\Rr^m
 \times [0,1]$ be defined as $\Phi(\omega)=(\omega, X(\omega),U(\omega))$ It is measurable for $\Fc$ and $\Ac \otimes \Rc^m \otimes \Rc_{[0,1]}$ where $\Rc^m$ is the Borel $\sigma-$algebra  on $\Rr^m$ and $\Rc_{[0,1]}$ the Borel $\sigma-$algebra on $[0,1]$.\\
 \item On $\Omega \times \Rr^m \times [0,1]$ we define $\chi(\omega,x,t)=x$ and $\tau(\omega,x,t)=t$. Clearly $\chi \circ \Phi=X$ and $\tau \circ \Phi=U.$\\
 \item The image probability $(\Pr\circ\Phi^{-1})$ of $\Phi$ can now be disintegrated, see \cite{Rao} for the existence of the conditional distributions. There exist kernels:\\
 $\mu_{X,U}:\Omega \times (\Rc^m\otimes \Rc_{[0,1]})\rightarrow[0,1],~\Ac-$measurable\\
 \\
 $\mu_X: \Omega\otimes \Rc^m\rightarrow[0,1],~ \Ac-$measurable, such that\\ \\
$\Pr\circ\Phi^{-1}(A \times B\times C)=\int_{A} \Pr[d\omega] \mu_{X,U}(\omega, B \times C), ~~ A\in \Ac, B\in\Rc^m, C\in \Rc_{[0,1]}$\\
 \\
 $\Pr\circ\Phi^{-1}(A \times B\times [0,1])=\int_{A} \Pr[d\omega] \mu_{X}(\omega,B),  A\in \Ac, B\in \Rc^m$
 \end{enumerate}
Almost surely the measure $\mu_X(\omega,.)$ is the marginal distribution of the probability measure $\mu_{X,U}(\omega,.)$.  Because of independence the conditional distribution of $U$ (given $\Ac$) is always the Lebesgue measure $\lambda$ on $[0,1]$.
\begin{theorem} We use the notation of the previous section. Let $Y\colon\Omega\rightarrow \Rr^m;Y\in L^2$ be independent of $\Ac$.  Let $\nu$ be the distribution of $Y$. Then
$$
\Vert X-Y\Vert^2\ge \int_\Omega \Pr[d\omega]\, W^2(\mu_X(\omega),\nu)
$$
\end{theorem}
 {\bf Proof}  Since $\Vert X\Vert^2=\Er[|X|^2]=\int_\Omega \Pr[d\omega]\int_{\Rr^m}|x|^2\mu_X(dx)$ we have that almost surely $\mu_X(\omega)\in\Mb_2(m)$. The disintegration of the distribution of $Y$ given $\Ac$ is --- because of independence --- trivially the constant $\nu$,  of course $\nu\in \Mb_2(m)$.
  The mapping $\Omega\rightarrow \Mb_2(m);\omega\rightarrow \mu_X(\omega)$ is $\Ac-$measurable when $\Mb_2(m)$ is endowed with its Borel $\sigma-$algebra. The integral in the statement of the theorem therefore makes sense since also $W(.,\nu)$ is Borel measurable, even continuous,  on $\Mb_2(m)$. Let $\mu_{X.Y}$ be the disintegration of the distribution of the couple $(X,Y)$ given $\Ac$. Since
  $$
  \Vert X-Y\Vert^2=\int_\Omega \Pr[d\omega]\int_{\Rr^m\times\Rr^m}|x-y|^2\mu_{X,Y}(\omega)(dx,dy),
  $$
  and since 
  $$
  \int_{\Rr^m\times\Rr^m}|x-y|^2\mu_{X,Y}(\omega)(dx,dy)\ge W^2(\mu_X(\omega),\nu)\text{ almost surely}
  $$
the theorem follows.\hfill$\Box$
\begin{theorem} There is a probability measure $\nu_0$ where the function $$\Mb_2(m)\rightarrow\Rr_+; \nu\rightarrow  \int_\Omega \Pr[d\omega]\, W^2(\mu_X(\omega),\nu)$$
attains its minimum.
\end{theorem}
{\bf Proof} Let $\nu_n$ be a sequence so that the integrals converge to the infimum.  The second moment of $\nu_n$ can be calculated using the Wasserstein distance to the Dirac measure  $\delta_0$, concentrated at $0\in\Rr^m$. Then we use $W(\nu_n,\delta_0)\le W(\mu_X,\delta_0)+W(\mu_X,\nu_n)$ and integrate. Since the sequence has uniformly bounded second moments, Prohorov's theorem implies that it is relatively weak$^*$ compact.  By selecting a subsequence we may suppose that it converges weak$^*$ to a probability measure $\nu_0$. Obviously we have that $\nu_0\in\Mb_2(m)$. The Wasserstein metric satisfies for any fixed $\alpha\in\Mb_2(m)$
$$
W(\alpha,\nu_0)\le \liminf W(\alpha,\nu_n).
$$
 From there it follows using Fatou's lemma that $ \int_\Omega \Pr[d\omega] W^2(\mu_X(\omega),\nu_0)$ realises the minimum.\hfill$\Box$

The only problem that remains, is to find a random variable that has distribution equal to $\nu_0$  that is independent of $\Ac$ and gives almost surely the Wasserstein distance to $\mu_X(\omega,.)$. Using the results from transport theory we will construct $Y$ for each $\omega\in\Omega$, or better for a set of full measure. $Y(\omega)$ will be a function of $X(\omega)$ and $U(\omega)$. This will guarantee measurability. $Y$ will have conditional to $\Ac$ and independent of $\omega$, the distribution $\nu_0$.  This will guarantee independence and it will in fact conclude the construction.  The details may look a little bit complicated but the construction is in fact the obvious way to follow.
 
\begin{theorem}  There is a random variable $Y\colon\Omega\rightarrow \Rr^m$ that is independent of $\Ac$, has $\nu_0$ as its distribution and realises the infimum 
$$
\Vert  X-Y\Vert=\inf\{ \Vert X-Z\Vert\mid Z\text{ independent of }\Ac\},
$$
\end{theorem}
{\bf Proof}  We use the functions $h$ and $\kappa$ from theorems \ref{h} and \ref{kappa} and make the reasoning for each $\omega$ separately. We start by looking at the measure $\mu_{X,U}$ on $\Rr^m\times[0,1]$.  On this set there are two obvious maps $\chi(x,t)=x$ having distribution $\mu_X(\omega)$ and $\tau(x,t)=t$ having distribution, the Lebesgue measure $\lambda$. We now use  on $\Rr^m\times\Rr^m$ the measure $\kappa(\mu_X(\omega),\nu_0)$.   We disintegrate this measure with respect to the first coordinate giving a kernel $\alpha(\omega,x,.)$ which is  measurable in $(\omega,x)$. To be precise we have
$$
\kappa(\mu_X(\omega),\nu_0)(A\times B)=\int_{A}\mu_X(\omega)(dx)\alpha(\omega,x,B).
$$
We now define $\upsilon(\omega,x,t)=h(\alpha(\omega,x),t)$. By construction and by definition of the kernel $\alpha$ we have that for each $\omega$ the map
$$
\upsilon(\omega)\colon \Rr^m\times [0,1]\rightarrow \Rr^m$$ 
has distribution $\nu_0$. The definition of $\kappa$ implies that
$$
W^2(\mu_X(\omega),\nu_0)=\int_{\Rr^m\times [0,1]} | \chi(x,t) - \upsilon(\omega,x,t) |^2\mu_{X,U}(dx,dt).
$$
We now put $Y=\upsilon\circ\Phi$ and recall that $X=\chi\circ \pi(\Phi)$ where $\pi$ maps $\Omega\times\Rr^m\times[0,1]$ onto $\Rr^m\times[0,1]$. The construction guarantees that $Y$ has the distribution $\nu_0$ and is independent of $\Ac$. As we saw this means that $Y$ is the best approximation of $X$ by a random variable independent of $\Ac$.\hfill$\Box$

 \section{A case with few independent random variables}
 Let us first recall the following definition.
\begin{definition} $\Fc$ is atomless conditionally to $\Ac$ if for every $A\in \Fc$, $\Pr[A]>0$, there is $B\subset A$ such that 
$$\Pr\[0< \Er[\one_B\mid\Ac]<\Er[\one_A\mid\Ac]\]>0.$$
\end{definition}
\begin{proposition} The existence of $U\colon\Omega\rightarrow \Rr$, uniformly distributed on $[0,1]$ and independent of $\Ac$, is equivalent to the condition that $\Fc$ is atomless conditionally to $\Ac$. \cite{del-atom} \end{proposition}

If no such condition exists then the existence of a uniformly $[0,1]$ distributed independent random variable, $U$, is no longer guaranteed. Nevertheless in some cases the existence of an optimal approximating random variable $Y$ can be shown to exist.  Let us illustrate this by the following rather extreme example.  Suppose $\Fc$ is generated by $\Ac$ and one set $A,0<\Pr[A]<1$. We suppose that $A$ is independent of $\Ac$ and look for other sets in $\Fc$ that are independent of $\Ac$. We first recall that
$$\Fc=\{(A \cap B)\cup(A^c\cap C)| B,C\in \Ac\}.$$
Suppose first that $\Pr[A]=\frac{1}{2}$ (from the general case it will follow that  this is indeed an isolated case) and take $B\in \Ac$ then
$$\Pr(A \cap B)\cup(A^c\cap B^c)=\Pr[A \cap B]+\Pr[A^c \cap B^c]= \frac{1}{2}\Pr[B]+\frac{1}{2}\Pr[B^c]=\frac{1}{2}$$
Now $D=(A \cap B)\cup(A^c \cap B^c)$ is independent of $\Ac$. Indeed
$$\Er[1_D\mid \Ac]=\Er[1_A\mid \Ac]1_B+\Er[1_A^c\mid \Ac]1_{B^c}=\frac{1}{2}$$
So the conditional expectation is constant and for sets this is equivalent to independence. All sets that are independent of $\Ac$ are of this form. Indeed for $(A \cap B)\cup(A^c \cap C)$ we have
$$\Er[1_D\mid \Ac]=\frac{1}{2}\one_B+\frac{1}{2}\one_C$$ 
and this is a constant only if $B=C=\Omega$ ($\Pr[D]=1$), $B=C=\emptyset$ ($\Pr[D]=0$) or $C=B^c$ and  then $\Pr[D]=\frac{1}{2}$.

Suppose that $\Dc$ is a sigma  algebra that is independent of $\Ac$ then $\Dc=\{\emptyset,\Omega, D=(A\cap  B)\cup(A^c\cap B^c), D^c\}$ where $B\in\Ac$. Indeed let $D_1=(A\cap  B_1)\cup(A^c\cap B_1^c), D_2=(A\cap  B_2)\cup(A^c\cap B_2^c)$ be two elements in $\Dc$, then unless we have trivial cases, $D_1\cap D_2$ is not of the form $(A\cap  C)\cup(A^c\cap C^c)$. It follows  that random variables $Y$ that are independent of  $\Ac$ can only take at most two different  values.  Suppose that $X$ is $\Fc$  measurable  then it can be written as  $X=f\one_A +g\one_{A^c}$ where $f,g$ are $\Ac$ measurable.  To avoid irrelevant complications we suppose that $f,g\ge 0$. Let $Y=\alpha \one_D +\beta\one_{D^c}$ where $D=(A\cap B)\cup(A^c\cap B^c)$  with $B\in\Ac$. The norm of the  difference $X-Y$ can now be calculated.
\begin{align*}
&\, \Vert X-Y\Vert^2\\
& = \Vert f\one_A +g\one_{A^c}-\alpha\(\one_{B\cap A}+\one_{b^c\cap A^c}\)-\beta\(\one_{B^c\cap A}+\one_{B\cap A^c} \)\Vert^2\\
& =  \Vert (f-\alpha)\one_{B\cap A}+(g-\alpha)\one_{B^c\cap A^c} + 
(f-\beta)\one_{B^c\cap A} +(g-\beta)\one_{B\cap A^c}\Vert^2\\
& =   \Er\[(f-\alpha)^2\one_{B\cap A}\]+\Er\[(g-\alpha)^2\one_{B^c\cap A^c}\] \\
& \quad + \,\Er\[(f-\beta)^2\one_{B^c\cap A}\] +\Er\[(g-\beta)^2\one_{B\cap A^c}\]\\
& =  \frac{1}{2}\(\Er\[(f-\alpha)^2\one_{B}\]+\Er\[(g-\alpha)^2\one_{B^c}\] \)\\
&\quad + \frac{1}{2} \( \Er\[(f-\beta)^2\one_{B^c}\] +\Er\[(g-\beta)^2\one_{B}\]\)\\
& =   \frac{1}{2}\(\Er\[(f-\alpha)^2\one_{B}+(g-\alpha)^2\one_{B^c}\] + 
\Er\[(f-\beta)^2\one_{B^c}+(g-\beta)^2\one_{B}\]\)\\
& =  \frac{1}{2}\(\Er\[(f\one_B+g\one_{B^c} - \alpha)^2\] + \Er\[(f\one_{B^c}+g\one_{B} - \beta)^2\]\).
\end{align*}
The optimal choices for $\alpha,\beta$ are
\begin{align*}
\alpha & = \Er[f\one_B+g\one_{B^c}]\\
\beta & = \Er[f\one_{B^c}+g\one_{B}],
\end{align*}
and then the expression becomes  (after rearrangement)
$$
\frac{1}{2}\(\Er[f^2]+\Er[g^2]- \(\Er[f\one_B+g\one_{B^c}]^2 +\Er[f\one_{B^c}+g\one_{B}]^2\)\).
$$
We can now choose the set $B$ so that the expression on the right becomes as big as possible.  An optimal choice is $B=\{f\ge g\}\in\Ac$ (because $f,g$ are nonnegative).  We then find
$$
\frac{1}{2}\(\Er[(f\vee g)^2]-\Er[(f\vee g)]^2 + \Er[(f\wedge g)^2]-\Er[(f\wedge g)]^2 \).
$$
For the optimal random  variable $Y$ we get:
$$
Y=
\Er[f\vee g]\one_A + \Er[f\wedge g]\one_{A^c}.
$$

Let us now quickly  analyse the  case $\Pr[A]\neq \frac{1}{2}$ then for
 $$D=(A \cap B)\cup(A^c \cap C)$$ 
 such that $D$ is independent of $\Ac$, we find
 $$\Pr[D]=\Pr[A]\Pr[B]+\Pr[A^c]\Pr[C]=\Pr[C]+\Pr[A](\Pr[B]-\Pr[C]).$$
 The conditional expectation is given by
 $$
 \Er[1_D\mid \Ac]=\Pr[A]\one_B+\Pr[A^c]\one_C
 $$
We first deal with the case $\Pr[D]=0$ which is only possible if $B=C=\emptyset$. Next we suppose that $\Pr[B\cap C]>0$. Then on $B\cap C$, the conditional expectation is equal to $1$ and since it must be a constant, it is $1$ on $\Omega$.  This means $\Pr[D]=1=\Pr[B]=\Pr[C]$. The remaining case is  when $\Pr[B\cap C]=0$. Then because  $\Pr[A]\neq\Pr[A^c]$ we must have that either $B$ or $C$ is the empty set and that $B=C^c$. That means that either $D=A$ or $D=A^c$. There is not much choice for a set (and hence for  a random variable) to be independent of $\Ac$.  The calculation of the optimal choice for $Y$ is left to the reader.


\begin{thebibliography}{100}

\bibitem{Black-Dub} Blackwell, D. and Dubins, L.E.: An Extension of Skorohod's almost sure Representation Theorem, {\em Proc. Amer. Math.Soc. 89 (4), 1983}

\bibitem{Bill} Billingsley, P.: Probability and Measure: Anniversary Edition, {\em Wiley, New Jersey, 2012}

\bibitem{del-atom} Delbaen, F.: Conditionally Atomless Extensions of Sigma Algebras, {\em arxiv 2003.09254, (2020) }

\bibitem{HLP} Hardy, G., Littlewood, J.E.,  Polya, G.: Inequalities {\em Cambridge Mathematical Library, Cambridge, 1934} 

\bibitem{Par} Parthasarathy, K.R.: Probability Measures on Metric Spaces, {\em Academic Press, New York, 1967}

\bibitem{Rao} Rao, M.M. and Swift, R.J.: Probability Theory with Applications, {\em 2nd Ed. Springer, New York, 2006}

\bibitem{Sri} Srivastava, S.M.: A Course on Borel Sets, {\em Springer, GTM 180, New York, 1998}

\bibitem{Thorpe} Thorpe, M.: Introduction to Optimal Transport {\em F2.08, Centre for Mathematical Sciences, University of Cambridge, 2018}

\end{thebibliography}
\end{document}